\documentclass[a4paper,12pt]{amsart}
\usepackage{amssymb}
\usepackage{ifthen}
\usepackage{graphicx}
\usepackage{subfigure}
\usepackage{color,xcolor}
\nonstopmode
\numberwithin{equation}{section}
\newtheorem{thm}{Theorem}
\newtheorem{cor}{Corollary}
\newtheorem{lem}{Lemma}
\newtheorem{prop}{Proposition}

\newtheorem{conj}{Conjecture}
\newtheorem{prob}{Problem}

\theoremstyle{definition}
\newtheorem{defn}{Definition}
\newtheorem{ca}{Case}

\newtheorem{rem}{Remark}

\newenvironment{pf}[1][]{%
 \vskip 1mm
 \noindent
 \ifthenelse{\equal{#1}{}}%
  {{\slshape Proof. }}%
  {{\slshape #1.} }%
 }%
{\qed\medskip}

\newcounter{alphabet}
\newcounter{tmp}
\newenvironment{Thm}[1][]{\refstepcounter{alphabet}%
\bigskip%
\noindent%
{\bf Theorem \Alph{alphabet}}%
\ifthenelse{\equal{#1}{}}{}{ (#1)}%
{\bf .} \itshape}{\vskip 8pt}

\makeatletter
\newcommand{\Ref}[1]{\@ifundefined{r@#1}{}{\setcounter{tmp}{\ref{#1}}\Alph{tmp}}}
\makeatother

\newcounter{alphabet2}

\newcommand{\IN}{{\mathbb N}}
\newcommand{\IC}{{\mathbb C}}
\newcommand{\ID}{{\mathbb D}}

\def\be{\begin{equation}}
\def\ee{\end{equation}}

\newcommand{\ben}{\begin{enumerate}}
\newcommand{\een}{\end{enumerate}}

\newcommand{\blem}{\begin{lem}}
\newcommand{\elem}{\end{lem}}
\newcommand{\bthm}{\begin{thm}}
\newcommand{\ethm}{\end{thm}}
\newcommand{\bcor}{\begin{cor}}
\newcommand{\ecor}{\end{cor}}
\newcommand{\beg}{\begin{exam}}
\newcommand{\eeg}{\end{exam}}
\newcommand{\begs}{\begin{examples}}
\newcommand{\eegs}{\end{examples}}
\newcommand{\bdefe}{\begin{defn}}
\newcommand{\edefe}{\end{defn}}
\newcommand{\bprob}{\begin{prob}}
\newcommand{\eprob}{\end{prob}}
\newcommand{\bques}{\begin{ques}}
\newcommand{\eques}{\end{ques}}
\newcommand{\bei}{\begin{itemize}}
\newcommand{\eei}{\end{itemize}}
\newcommand{\bcon}{\begin{conj}}
\newcommand{\econ}{\end{conj}}
\newcommand{\bop}{\begin{op}}
\newcommand{\eop}{\end{op}}

\newcommand{\bas}{\begin{assertion}}
\newcommand{\eas}{\end{assertion}}

\newcommand{\bfa}{\begin{fact}}
\newcommand{\efa}{\end{fact}}

\newcommand{\bca}{\begin{ca}}
\newcommand{\eca}{\end{ca}}

\newcommand{\bst}{\begin{step}}
\newcommand{\est}{\end{step}}

\newcommand{\bsca}{\begin{sca}}
\newcommand{\esca}{\end{sca}}

\newcommand{\bcl}{\begin{cl}}
\newcommand{\ecl}{\end{cl}}

\newcommand{\bmlem}{\begin{mlem}}
\newcommand{\emlem}{\end{mlem}}

\newcommand{\bscl}{\begin{scl}}
\newcommand{\escl}{\end{scl}}

\newcommand{\bcons}{\begin{conjs}}
\newcommand{\econs}{\end{conjs}}

\newcommand{\bprop}{\begin{prop}}
\newcommand{\eprop}{\end{prop}}

\newcommand{\br}{\begin{rem}}
\newcommand{\er}{\end{rem}}
\newcommand{\brs}{\begin{rems}}
\newcommand{\ers}{\end{rems}}
\newcommand{\bo}{\begin{obser}}
\newcommand{\eo}{\end{obser}}
\newcommand{\bos}{\begin{obsers}}
\newcommand{\eos}{\end{obsers}}
\newcommand{\bpf}{\begin{pf}}
\newcommand{\epf}{\end{pf}}
\newcommand{\ba}{\begin{array}}
\newcommand{\ea}{\end{array}}
\newcommand{\beq}{\begin{eqnarray}}
\newcommand{\beqq}{\begin{eqnarray*}}
\newcommand{\eeq}{\end{eqnarray}}
\newcommand{\eeqq}{\end{eqnarray*}}

\newcommand{\ds}{\displaystyle}

\newcounter{minutes}
\divide\time by 60
\newcounter{hours}
\multiply\time by 60 \addtocounter{minutes}{-\time}

\begin{document}

\bibliographystyle{amsplain}
\title [Refined bohr-type inequalities for bounded analytic functions]
{Refined bohr-type inequalities with area measure for bounded analytic functions}

\def\thefootnote{}
\footnotetext{ \texttt{\tiny File:~\jobname .tex,
          printed: \number\day-\number\month-\number\year,
          \thehours.\ifnum\theminutes<10{0}\fi\theminutes}
} \makeatletter\def\thefootnote{\@arabic\c@footnote}\makeatother

\author{Yong Huang}
 \address{Y. Huang, School of Mathematical Sciences, South China Normal University, Guangzhou, Guangdong 510631, China.}
 \email{hyong95@163.com}

\author{Ming-Sheng Liu${}^{~\mathbf{*}}$}
 \address{M-S Liu, School of Mathematical Sciences, South China Normal University, Guangzhou, Guangdong 510631, China.} \email{liumsh65@163.com}

\author{Saminathan Ponnusamy}
\address{S. Ponnusamy, Department of Mathematics,
Indian Institute of Technology Madras, Chennai-600 036, India. }
\email{samy@iitm.ac.in}



\subjclass[2000]{Primary: 30A10, 30C45, 30C62; Secondary: 30C75}
\keywords{Bohr radius, bounded analytic functions, harmonic function, Bohr inequality
  \\
${}^{\mathbf{*}}$ Correspondence should be addressed to Ming-Sheng Liu
}


\begin{abstract}
In this paper, we establish five new sharp versions of Bohr-type inequalities for bounded analytic functions in the unit disk by allowing  Schwarz function
in place of the initial coefficients in the power series representations of the functions involved and thereby, we generalize several related results of earlier authors.
\end{abstract}

\maketitle
\pagestyle{myheadings}
\markboth{Y. Huang, M-S Liu and S. Ponnusamy}{Refined bohr-type inequalities for bounded analytic functions}

\section{Introduction and Preliminaries}\label{HLP-sec1}

Let $\mathbb{D}:=\{z \in \mathbb{C}:\,|z|<1\}$ denote the open unit disk in $\IC$. A remarkable discovery of Herald Bohr \cite{B1914} in 1914
states that if $H_\infty $ denotes the class of all bounded analytic functions $f$ on $\ID$ with the supremum norm
$\|f\|_\infty :=\sup_{z\in \ID}|f(z)|$, then
\begin{equation}
B_0(f,r):= |a_0|+\sum_{n=1}^{\infty} |a_n| r^n\leq \|f\|_\infty ~\mbox{ for $0\leq r\leq 1/6$,}
\label{liu1}
\end{equation}
where $a_k=f^{(k)}(0)/k!$ for $k\geq 0$. Later M.~Riesz, I.~Shur and F. W. Wiener, independently proved its validity on a wider
interval $0\leq r\leq 1/3$, and 
the family of functions $\varphi_a(z)=(a-z)/(1-\overline{a} z)$ ($|a|<1$) as $a\rightarrow 1$ demonstrates that the number $1/3$ is optimal.
This result is usually referred to as Bohr's power series theorem for the unit
disk and $1/3$ is called the Bohr radius. We refer the paper of Bohr \cite{B1914} which contains the proof of Wiener showing that the
Bohr radius is $1/3$. See also \cite{S1927,T1962} for other proofs. Then it is worth pointing out that there is no extremal function in $\mathcal{B}$
such that the Bohr radius is precisely $1/3$ (cf. \cite{AlKayPON-19}, \cite[Corollary 8.26]{GarMasRoss-2018} and \cite{KP2017}). 
Several aspects  of Bohr's inequality and its extensions in various settings may be seen in the literature.
For example, the Bohr radius for analytic functions from the unit disk into different domains, such as the punctured unit disk or the exterior of the closed unit disk
or concave wedge-domains,  have been analyzied in \cite{Abu,Abu2,Abu4,Abu3}.
Ali et al. \cite{AliBarSoly, KP2017} considered the problem of determining Bohr radius for the classes of even and odd analytic functions and for alternating series.
The articles \cite{AlKayPON-19,KSS2017,LP2019} concerned with the class of all sense-preserving harmonic mappings and the Bohr radius for  sense-preserving harmonic
quasiconformal mappings. Defant \cite{DFOOS} improved a version of the Bohnenblust-Hille inequality, and in 2004, Paulsen \cite{PS2004} proved a uniform algebra analogue
of the classical inequality of Bohr concerning Fourier coefficients of bounded holomorphic functions. In \cite{PPS2002,PS2006}, the authors demonstrated
the classical Bohr inequality using different methods of operators. Djakov and Ramanujan \cite{DjaRaman-2000} have established the results on Bohr's phenomena for multidimensional
power series. Recently, in \cite{LSX2018,LPW2020,PVW2019,PW2019}, the authors presented refined versions of Bohr's inequality along with few other related
improved versions of previously known results. See also the recent survey articles \cite{AAP2016,IKKP2018,KKP2018} and \cite[Chapter 8]{GarMasRoss-2018}.
Especially, after the appearance of the articles \cite{AAP2016} and  \cite{KS2017},
several approaches and new problems on Bohr's inequality in the plane were investigated in the literature
(cf. \cite{BenDahKha,BhowDas-18,KayPon3,LP2018, LSX2018,PVW2019,PW2019}).

One of our aims in this article is to generalize or improve recent versions of Bohr's inequalities for functions from $H_\infty $.

\subsection{Basic Notations}
Before we continue the discussion, we fix some notations. Throughout the discussion, we let
\beqq
{\mathcal B}& = &\{f\in H_\infty :\, \|f\|_\infty \leq 1 \}, \mbox{ and }  ~m\in \mathbb{N}=\{1,2,\cdots\}, \\
{\mathcal B}_m&=&\{\omega \in {\mathcal B}:\, \omega (0)= \cdots =\omega ^{(m-1)}(0)=0 ~\mbox{ and }~ \omega ^{(m)}(0)\neq 0 \}.
\eeqq
Also, for $f(z)=\sum_{n=0}^{\infty} a_{n} z^{n}\in {\mathcal B}$ and $f_0(z):=f(z)-f(0)$, we let for convenience
$$B_{k}(f,r) :=  \sum_{n=k}^{\infty} |a_n| r^n ~~\mbox{for $k\geq 0$,} ~\mbox{ and }~
\|f_0\|_r^2  :=   \sum_{n=1}^{\infty}\left|a_{n}\right|^{2} r^{2 n}
$$
so that $B_0(f,r) = |a_0|+ B_1(f,r)$ and $B_0(f,r) = |a_0|+ |a_1|r+B_2(f,r)$.

\subsection{Refined Bohr's inequalities and basic problems}
Recently, Ponnusamy  et al. \cite{PVW2019} proved the following refined Bohr inequality.

\begin{Thm}\label{Theo-B}
(\cite{PVW2019})
Suppose that $f \in\mathcal{B}$, $f(z)=\sum_{n=0}^{\infty} a_{n} z^{n}$ and $a=|a_0|=|f(0)|$. Then
$$
B_0(f,r)+ \frac{1+ar}{(1+a)(1-r)}\|f_0\|_r^2  \leq 1 ~\mbox{ for }~  r \leq \frac{1}{2 +a}
$$
and the numbers $ \frac{1}{2+a}$ and  $\frac{1}{1+a}$ cannot be improved. Moreover,
$$
 a ^{2}+B_1(f,r)+\frac{1+ar}{(1+a)(1-r)}\|f_0\|_r^2 \leq 1  ~for~ r \leq \frac{1}{2}
$$
and the numbers $\frac{1}{2}$ and  $\frac{1}{1+ a}$ cannot be improved.
\end{Thm}

Besides these results, there are plenty of works about the classical Bohr inequality.
Based on the work of Kayumov and Ponnusamy \cite{KayPon3}, several forms of Bohr-type inequalities for the family $\mathcal{B}$ were
obtained in \cite{LSX2018} when the Taylor coefficients of classical Bohr inequality are partly or completely replaced by higher order derivatives of $f$.
With the development of Bohr-type inequalities, the authors in \cite{LLP2020}
established improved version of the Bohr-Rogosinski inequality  and considered some refined Bohr type inequalities
associated with area, modulus of $f-a_{0}(f)$ and higher order derivatives of $f$ in part. Here we recall part of them.

\begin{Thm}\label{Theo-C}
(\cite{LLP2020})
Suppose that $f \in\mathcal{B}$, $f(z)=\sum_{n=0}^{\infty} a_{n} z^{n}$ and $a:=|a_0|=|f(0)|$. Then
$$
|f(z)|+B_1(f,r)+ \frac{1+ar}{(1+a)(1-r)} \|f_0\|_r^2 \leq 1
$$
for $|z|=r \leq r_{a}=2 /\left(3+a+\sqrt{5}\left(1+a\right)\right)$. Then the radius $ r_{a}$ is best possible and $r_{a} \geq \sqrt{5}-2$.
Moreover,
$$
|f(z)|^2+B_1(f,r)+ \frac{1+ar}{(1+a)(1-r)} \|f_0\|_r^2 \leq 1
$$
for $|z|=r \leq r_{a}^{\prime}$, where $r_{a}^{\prime}$ is the unique positive root of the equation
$$
\left(1-a^{3}\right) r^{3}-\left(1+2a\right) r^{2}-2 r+1=0.
$$
The radius $r_{a}^{\prime}$ is best possible. Further,we have $1/3 <r_{a}^{\prime}<1 /\left(2+\left|a\right|\right)$.
\end{Thm}

\begin{Thm}\label{Theo-D}
(\cite{LLP2020})
Suppose that $f \in\mathcal{B}$, $f(z)=\sum_{n=0}^{\infty} a_{n} z^{n}$ and $a:=|a_0|=|f(0)|$. Then
$$
B_0(f,r)+ \frac{1+ar}{(1+a)(1-r)} \|f_0\|_r^2+ |f_0(z)| \leq 1
$$
for $|z|=r \leq \frac{1}{5}$ and the number $\frac{1}{5}$ cannot be improved. Moreover,
$$a^{2}+B_1(f,r)+ \frac{1+ar}{(1+a)(1-r)} \|f_0\|_r^2+|f_0(z)| \leq 1  ~
$$
for $|z|=r \leq \frac{1}{3}$ and the constant $\frac{1}{3}$ cannot be improved.
\end{Thm}

\begin{Thm}\label{Theo-E}
(\cite{LLP2020})
Suppose that $f \in\mathcal{B}$, $f(z)=\sum_{n=0}^{\infty} a_{n} z^{n}$ and $a:=|a_0|=|f(0)|$. Then
$$
|f(z)|+\left|f^{\prime}(z)\right| r+B_2(f,r)+ \frac{1+ar}{(1+a)(1-r)} \|f_0\|_r^2 \leq 1
$$
for $|z|=r \leq \frac{\sqrt{17}-3}{4}$ and the constant $\frac{\sqrt{17}-3}{4}$ is best possible.Moreover,
$$
|f(z)|^2+\left|f^{\prime}(z)\right| r+B_2(f,r)+ \frac{1+ar}{(1+a)(1-r)} \|f_0\|_r^2 \leq 1
$$
for $|z|=r \leq r_{0}$, where $ r_{0}\approx 0.385795$ is the unique positive root of the equation
$$
1-2r-r^{2}-r^{3}-r^{4}=0
$$
and the number $ r_{0}$ is best possible.
\end{Thm}

\begin{Thm}\label{Theo-F}
(\cite{LLP2020})
Suppose that $f \in\mathcal{B}$, $f(z)=\sum_{n=0}^{\infty} a_{n} z^{n}$ and $a:=|a_0|=|f(0)|$. Then for $|z|=r\leq \frac{1}{3}$,
$$
B_0(f,r)+ \frac{1+ar}{(1+a)(1-r)} \|f_0\|_r^2+|f_0(z)|^{2} \leq 1 ~
$$
is valid if and only if $0\leq a\leq 4\sqrt{2}-5\approx 0.656854$.
\end{Thm}

It is natural to raise the following.

\bprob\label{HLP-prob1}
Whether we can further generalize or improve Theorems \Ref{Theo-C}, \Ref{Theo-D}, \Ref{Theo-E} and \Ref{Theo-F}?
\eprob

In this article, we present an affirmative answer to this question in five different forms.

The paper is organized as follows. In Section \ref{HLP-sec2}, we present statements of five theorems which improve several versions of
Bohr's type inequalities for bounded analytic functions, and several remarks. In Section \ref{HLP-sec3}, we state and prove a couple of lemmas
which are needed for the proofs of two theorems.
In Section \ref{HLP-sec4}, we present the proofs of the main results.

\section{Statement of Main Results and Remarks}\label{HLP-sec2}
We now state a generalization of Theorem \Ref{Theo-C} in a general setting.

\bthm\label{HLP-th1}
Suppose that $f(z)=\sum_{n=0}^{\infty} a_{n} z^{n} \in\mathcal{B}$,  $a:=|a_0|$ and $\omega \in\mathcal{B}_m$ for some $m \in \IN$.
Then we have
\begin{eqnarray*}
A_f(z):=\left|f\left(\omega(z)\right)\right|+B_1(f,r)+ \frac{1+ar}{(1+a)(1-r)} \|f_0\|_r^2 \leq 1
\end{eqnarray*}
for $r \in\left[0, \alpha_{m}\right]$, where $\alpha_{m}$ is the unique root in $(0, 1)$ of the equation
\begin{equation}
(1-r)(1-r^m)-2r(1+r^m)=0.
\label{liu21}
\end{equation}
The constant $\alpha_{m}$ cannot be improved. Moreover,
\begin{eqnarray*}
B_f(z):=\left|f\left(\omega(z)\right)\right|^2+B_1(f,r)+ \frac{1+ar}{(1+a)(1-r)} \|f_0\|_r^2 \leq 1
\end{eqnarray*}
is valid for $r \in\left[0, \beta_{m}\right] $, where $\beta_{m}$ is the unique root in $(0, 1)$ of the equation
\begin{equation}
1-2r-r^{m}=0 . 
\label{liu22}
\end{equation}
The constant $\beta_{m}$ cannot be improved.
\ethm

\br
We mention now several useful remarks and some special cases.
\begin{enumerate}
\item One can state each of the two radii in Theorem \ref{HLP-th1} as a function of $a$. In that case, $\alpha_{m}$ and $\beta_{m}$ should be
replaced by $\alpha_{m, a}$ and $\beta_{m,a}$ which are in fact
the unique roots in $(0, 1)$ of the equation $A_m(a,r)=0$ and $B_m(a,r)=0$, respectively, where
$$ A_m(a,r)=(1-r)(1-r^m)-(1+a)r(1+ar^m)
$$
and
$$
B_{m}(a, r)=(1-r)(1-r^{2m})-r(1+ar^m)^2.
$$
\item If we set $m=1$ and $\omega(z)\equiv z$ in Theorem \ref{HLP-th1}, then we get Theorem \Ref{Theo-C}.
\item  If we set $m=1$ in \eqref{liu21}, then we get  $\alpha_{1}=\sqrt{5} -2$.
\item  If we set $m=2$ in \eqref{liu21}, then it reduces to $ r^{3}+r^2+3 r -1=0$ which gives the root $\alpha_{2}\approx 0.295598$ in the interval $(0, 1)$.
\item  If we allow $m\rightarrow \infty $  in \eqref{liu21}  (with $\omega(z)=z^m$ in $A_f(z)$), then $|f(\omega(z))|\rightarrow |f(0)|$ and
$\alpha_{\infty}=1/3$.
\item  The case $m=1$ in \eqref{liu22} gives the root $\beta_{1}=1/3 $.
\item The case $m=2$ in \eqref{liu22}  gives the root $\beta_{2} =\sqrt{2} -1$.
\item If we allow $m\rightarrow \infty $ in \eqref{liu22} (with $\omega(z)=z^m$ in $B_f(z)$), then $\beta_\infty =1/2$.
\end{enumerate}
\er

\bthm \label{HLP-th2}
Suppose that $f(z)=\sum_{n=0}^{\infty} a_{n} z^{n} \in\mathcal{B}$,  $a:=|a_0|$ and $\omega \in\mathcal{B}_m$ for some $m \in \IN$.
Then we have
\begin{eqnarray*}
C_f(z):=B_0(f,r)+ \frac{1+ar}{(1+a)(1-r)} \|f_0\|_r^2 + \left|f\left(\omega(z)\right)-a_{0}\right| \leq 1
\end{eqnarray*}
for $r \in\left[0, \zeta_{m}\right] $, where $\zeta_{m}$ is the unique root in $(0, 1/3]$ of the equation
\begin{equation}
 r^{m}(3-5r)+3 r -1=0,
\label{liu23}
\end{equation}
or equivalently, $3r^m +2\sum_{k=1}^{m}r^{k}-1=0.$ The upper bound $\zeta_{m}$ cannot be improved.

Moreover,
\begin{eqnarray*}
D_f(z):=|a_{0}|^2+B_1(f,r)+ \frac{1+ar}{(1+a)(1-r)} \|f_0\|_r^2+ |f(\omega(z))-a_{0}| \leq 1
\end{eqnarray*}
for $r \in\left[0, \eta_{m}\right] $, where $\eta_{m}$ is the unique root in $(0, 1/2]$ of the equation
\begin{equation}
r^{m}(2-3r) +2 r-1=0, 
\label{liu24}
\end{equation}\
or equivalently, $2r^m +\sum_{k=1}^{m}r^{k}-1=0.$
The upper bound $\eta_{m}$ cannot be improved.
\ethm

\br
The following special cases are useful and important to mention.
\begin{enumerate}
\item The case $m=1$ and $\omega(z)\equiv z$ in Theorem \ref{HLP-th2} gives Theorem \Ref{Theo-D}.
\item  The case $m=1$ in \eqref{liu23} gives the root $\zeta_{1}=1/5$.
\item The case $m=2$ in \eqref{liu23} reduces to $ -5r^{3}+3r^2+3 r -1=(1-r)(5r^{2}+2r -1)=0$ which gives the root
$\zeta_{2}=\frac{\sqrt{6}-1}{5} \approx 0.289898 $ in the interval $(0, 1/3)$.
\item If we allow $m\rightarrow \infty $ in \eqref{liu23} (with $\omega(z)=z^m$ in $C_f(z)$), then $\zeta_\infty =1/3$.
\item  The case $m=1$ in \eqref{liu24} gives the root $\eta_{1}=1/3$.
\item The case $m=2$ in \eqref{liu24} reduces to $ (1-r)(3r^{2}+ r -1)=0$ which gives the root $\eta_{2} =\frac{\sqrt{13}-1}{6} \approx 0.434259$ in the interval $(0, 1/2)$.
\item If we allow $m\rightarrow \infty $ in \eqref{liu24} (with $\omega(z)=z^m$ in $D_f(z)$), then $\eta_\infty =1/2$.
\end{enumerate}

In Table \ref{Tab1ForThe1-2}, we include the values of $\alpha_{m}$, $\beta_{m}$, $\zeta_{m}$ and $\eta_{m}$ for certain values of $m\geq 3$.
\er

\begin{table}[tbp]
\centering
\begin{tabular}{|l||l|l||l|l|l|}
\hline
$m$  &$\alpha_{m}$  &$\beta_{m}$  &$\zeta_{m}$  &$\eta_{m}$ \\
\hline
3  &0.319053  &0.453398  &0.318201  &0.469396   \\
\hline
4  &0.328197  &0.474627  &0.328083  &0.484925    \\
\hline
5  &0.331555  &0.486389  &0.331541  &0.492432    \\
\hline
6  &0.332731  &0.492836  &0.332729  &0.496184    \\
\hline
7  &0.333131  &0.496292  &0.333131  &0.498077    \\
\hline
8  &0.333266  &0.498105  &0.333266  &0.499033    \\
\hline
9  &0.333311  &0.499040  &0.333311  &0.499515    \\
\hline
10 &0.333326  &0.499516  &0.333326  &0.499757    \\
\hline
15 &0.333333  &0.499985  &0.333333  &0.499992    \\
\hline
20 &0.333333  &0.500000  &0.333333  &0.500000    \\
\hline
25 &0.333333  &0.500000  &0.333333  &0.500000    \\
\hline
30 &0.333333  &0.500000  &0.333333  &0.500000    \\
\hline
\end{tabular}
\vspace{8pt}
\caption{Numbers $\alpha_{m}$, $\beta_{m}$, $\zeta_{m}$ and $\eta_{m}$ are the unique roots in $(0,1)$ of the equations
\eqref{liu21}, \eqref{liu22}, \eqref{liu23} and \eqref{liu24}, respectively. \label{Tab1ForThe1-2}
}
\end{table}

\bthm\label{HLP-th3}
Suppose that $f(z)=\sum_{n=0}^{\infty} a_{n} z^{n} \in\mathcal{B}$,  $a:=|a_0|$ and $\omega \in\mathcal{B}_m$.  Then we have
\begin{eqnarray*}
E_f(z):=\left|f\left(\omega(z)\right)\right|+|\omega(z) |\, |f^{\prime}(\omega(z))|+B_2(f,r)+\frac{1+ar}{(1+a)(1-r)} \|f_0\|_r^2 \leq 1
\end{eqnarray*}
for $r \in\left[0, \gamma_{m}\right] $, where $\gamma_{m}$ is the unique root in $(0, 1)$ of the equation
\begin{equation}
r^{m}(r^m+2)[2r^{2} -r +1]+2 r^{2}+r-1=0.
\label{liu25}
\end{equation}
The upper bound $\gamma_{m}$ cannot be improved. Moreover,
\begin{eqnarray*}
F_f(z):=\left|f\left(\omega(z)\right)\right|^2+|\omega(z) |\,| f^{\prime}(\omega(z))|+B_2(f,r)+\frac{1+ar}{(1+a)(1-r)} \|f_0\|_r^2 \leq 1
\end{eqnarray*}
for $r \in\left[0, \delta_{m}\right] $, where $\delta_{m}$ is the unique root in $(0, 1)$ of the equation
\begin{equation}
r^{m}(r^m+1)[r^{m} -r +2]+ r-1 =0.
\label{liu26}
\end{equation}
The upper bound $\delta_{m}$ cannot be improved.
\ethm

\bthm\label{HLP-th4}
Suppose that $f(z)=\sum_{n=0}^{\infty} a_{n} z^{n} \in\mathcal{B}$,  $a:=|a_0|$ and $\omega \in\mathcal{B}_m$.  Then we have
\begin{eqnarray*}
G_f(z):=\left|f\left(\omega(z)\right)\right|+|z |\,|f^{\prime}(\omega(z))|+B_2(f,r)+\frac{1+ar}{(1+a)(1-r)} \|f_0\|_r^2 \leq 1
\end{eqnarray*}
for $r \in\left[0, \theta_{m}\right] $, where $\theta_{m}$ is the unique root in $(0, 1)$ of the equation
\begin{equation}
2 r^{2 m+2}-r^{2 m+1}+r^{2 m}+4 r^{m+2}+3r-1=0.
\label{liu27}
\end{equation}
The upper bound $\theta_{m}$ cannot be improved. Moreover,
\begin{eqnarray*}
H_f(z):=\left|f\left(\omega(z)\right)\right|^2+ |z | \,| f^{\prime} (\omega(z))|+B_2(f,r)+\frac{1+ar}{(1+a)(1-r)} \|f_0\|_r^2 \leq 1
\end{eqnarray*}
for $r \in\left[0, \vartheta_{m}\right] $, where $\vartheta_{m}$ is the unique root in $(0, 1)$ of the equation
\begin{equation}
r^{2 m+2}-r^{2 m+1}+r^{2 m}+2 r^{m+2}+2r-1=0.
\label{liu28}
\end{equation}
The upper bound $\vartheta_{m}$ cannot be improved.
\ethm

\br
Obviously, if we set $m=1$ and $\omega(z)\equiv z$ in Theorem \ref{HLP-th3} or Theorem \ref{HLP-th4}, then we get Theorem \Ref{Theo-E}.

In Table \ref{Tab2ForThe3-4}, we include the values of $\gamma_{m}$, $\delta_{m}$, $\theta_{m}$ and $\vartheta_{m}$ for certain values of $m\geq 2$.
\er

\begin{table}[tbp]
\centering
\begin{tabular}{|l||l|l||l|l|}
\hline
$m$  &$\gamma_{m}$  &$\delta_{m}$  &$\theta_{m}$  &$\vartheta_{m}$\\
\hline %
2    &0.391490  &0.486848  &0.316912  &0.445688  \\
\hline
3    &0.441112  &0.535687  &0.327911  &0.472325  \\
\hline
4    &0.467644  &0.564540 &0.331520  &0.485708   \\
\hline
5    &0.482442  &0.582935  &0.332726  &0.492642  \\
\hline
6    &0.490660  &0.595034  &0.333131  &0.496239  \\
\hline
7    &0.495127  &0.603062  &0.333266  &0.498091  \\
\hline
8    &0.497496  &0.608373  &0.333311  &0.499037  \\
\hline
9    &0.498727  &0.611827  &0.333326  &0.499515  \\
\hline
10   &0.499357  &0.614117  &0.333331  &0.499757  \\
\hline
15   &0.499980  &0.617662  &0.333333  &0.499992  \\
\hline
20   &0.500000  &0.618000  &0.333333  &0.500000  \\
\hline
25   &0.500000  &0.618031  &0.333333  &0.500000  \\
\hline
30   &0.500000  &0.618034  &0.333333  &0.500000  \\
\hline
\end{tabular}
\vspace{8pt}
\caption{Numbers $\gamma_{m}$, $\delta_{m}$, $\theta_{m}$ and $\vartheta_{m}$ are the unique roots in $(0, 1)$ of the equations
\eqref{liu25}, \eqref{liu26}, \eqref{liu27} and \eqref{liu28}, respectively. \label{Tab2ForThe3-4}
}

\end{table}

\bthm \label{HLP-th5}
Suppose that $f(z)=\sum_{n=0}^{\infty} a_{n} z^{n} \in\mathcal{B}$,  $a:=|a_0|$ and $\omega \in\mathcal{B}_m$ for some $m\geq 1$.  We have the following:
\begin{enumerate}
\item[{\rm (1)}] If $m=1$, then we have
\be
 I_f(z):=B_0(f,r)+\frac{1+ar}{(1+a)(1-r)} \|f_0\|_r^2+\left|f\left(\omega(z)\right)-a_{0}\right|^{2} \leq 1 ~\mbox{ for $ |z|=r\leq 1/3 $}
\label{liu29}
\ee
if and only if $0\leq a \leq a^{*}=-5 +4\sqrt{2}\approx 0.656854$. The constant $1/3$ cannot be improved.

\item[{\rm (2)}] If $m\geq 2$, then \eqref{liu29} holds, and the constant $1/3$ cannot be improved.
\end{enumerate}
\ethm

\br
Obviously, the case $m=1$ and $\omega(z)\equiv z$ of Theorem \ref{HLP-th5} gives  Theorem \Ref{Theo-F}.
\er

\section{Key lemmas and their Proofs}\label{HLP-sec3}

In order to establish our main results, we need the several lemmas which play key role in proving the subsequent
results in Section \ref{HLP-sec4}.

\subsection{Three known lemmas}

\blem\label{HLP-lem5}(Schwarz-Pick Lemma)
Let $\varphi(z)$ be analytic and $|\varphi(z)|< 1$ in the unit disk $\mathbb{D}$. Then
$$\frac{\left|\varphi(z_1)-\varphi(z_2)\right|}{\left|1-\overline{\varphi(z_1)}\varphi(z_2)\right|}\leq\frac{\left|z_1-z_2\right|}{\left|1-\overline{z_1}z_2\right|}
~\mbox{ for $z_1, z_2\in \mathbb{D}$},
$$
and equality holds for distinct $z_1, z_2\in \mathbb{D}$ if and only if $\varphi$ is a {\it M$\ddot{o}$bius transformation}. Also,
$$|\varphi '(z)|\leq \frac{1-|\varphi(z)|^2}{1-|z|^2} ~\mbox{ for $z\in \mathbb{D}$},
$$
and equality holds for some $z\in \mathbb{D}$ if and only if $f$ is a {\it     M$\ddot{o}$bius transformation}.
\elem

\blem\label{HLP-lem1} (\cite{KP2017})
Suppose that $f(z)=\sum_{n=0}^{\infty} a_{n} z^{n} \in\mathcal{B}$ and  $a:=|a_0|$. Then we have
\begin{equation*}
\sum_{n=1}^{\infty}\left|a_{n}\right| r^{n} \leq\left \{ \begin{array}{lr} \ds r\frac{1-a^2}{1-ra} & \mbox{ for $a \ge r$}, \\[4mm]
\ds r\frac{\sqrt{1-a^2}}{\sqrt{1-r^{2}}} &\mbox{ for $a < r$}.
\end{array} \right .
\end{equation*}
\elem

\blem\label{HLP-lem2} (\cite{PVW2019})
Suppose that $f(z)=\sum_{n=0}^{\infty} a_{n} z^{n} \in\mathcal{B}$ and  $a:=|a_0|$ Then the following inequality holds:
$$ B_1(f,r)+ \frac{1+ar}{(1+a)(1-r)} \|f_0\|_r^2  \leq \left(1-a^{2}\right) \frac{r}{1-r} ~\mbox{ for $r\in[0,1)$.}~
$$
\elem

A general version of this lemma is proved in \cite[Lemma 4]{LLP2020}. In particular, the following inequality holds  (the case
$N=2$ in \cite[Lemma 4]{LLP2020}) under the hypothesis of Lemma \ref{HLP-lem2}:
\be\label{HLP-eq2}
B_2(f,r)+  \frac{1+ar}{(1+a)(1-r)} \|f_0\|_r^2 \leq\left(1-a^{2}\right) \frac{r^{2}}{1-r} ~\mbox{ for $r\in[0,1)$.}
\ee

%

\subsection{Two key lemmas}

\blem\label{HLP-lem3}
There is a unique positive root $\zeta_{m}$ in $(0, 1/3)$ of the equation \eqref{liu23}, and $\zeta_{m}$ satisfies the inequality
\begin{eqnarray}
\zeta_{m}^{m}+\frac{\zeta_{m}}{1-\zeta_{m}}+\frac{\zeta_{m}^{m}}{\sqrt{1-\zeta_{m}^{2 m}}}\leq 1.
\label{liu31}
\end{eqnarray}
\elem
\bpf
We first prove the uniqueness of the solution in $(0, 1/3)$ of the equation (\ref{liu23}).

Let $g(r)= r^{m}(3-5r)+3 r -1$. Then, we find that $g(0)= - 1<0$ and $g(1/3)=\frac{4}{3^{m+1}}> 0$. Also, for $r\in [0, 1/3]$,
we have
\begin{eqnarray*}
g'(r)= m r^{m-1} (3-5r)+3-5r^{m}>0,
\end{eqnarray*}
showing that $g(r)$ is an increasing function of $r$ in $[0, 1/3]$, and thus, $g(r)=0$ has a unique root $\zeta_{m}$ in $(0, 1/3)$.

Now we verify the inequality (\ref{liu31}). In fact, by (\ref{liu23}), we note that
$$\zeta_{m}^{m}=\frac{1-3 \zeta_{m}}{3-5 \zeta_{m}},\, \zeta_{m}\in(0, 1/3).
$$
For convenience, we set $x=\zeta_{m}$ and using the last relation, we have
\begin{eqnarray*}
x^{m}+\frac{x}{1-x}+\frac{x^{m}}{\sqrt{1-x^{2 m}}}&=&\frac{1-3 x}{3-5 x}\left(1+\frac{1}{\sqrt{1-\left(\frac{1-3 x}{3-5 x}\right)^{2}}}\right)+\frac{x}{1-x}\\
&\leq& \frac{1-3 x}{3-5 x}\left(1+\frac{1}{\sqrt{1-\left(\frac{1}{3}\right)^{2}}}\right)+\frac{x}{1-x}\\
&=&\frac{4+3 \sqrt{2}}{4} \cdot \frac{1-3 x}{3-5 x}+\frac{x}{1-x}\\
&=&\frac{(9\sqrt{2}-8) x^{2}-(12 \sqrt{2}+4) x+4+3 \sqrt{2}}{4\left(3-5 x\right)\left(1-x\right)}  ,
\end{eqnarray*}
which is less than or equal to $1$ if $G(x)\geq 0$, where
$$
G(x)=(28-9\sqrt{2}) x^{2}+(12 \sqrt{2}-28) x+8-3 \sqrt{2}.
$$

Since the discriminant of the equation $G(x)=0$ is less than 0 and $G(0)=8-3 \sqrt{2}>0$,
we deduce that $G(x)> 0$ for $x\in(0, 1/3)$. The proof is complete.
\epf

\blem\label{HLP-lem4}
There is a unique positive root $\eta_{m}$ in $(0, 1/2)$ of the equation \eqref{liu24}, and $\eta_{m}$ satisfies the inequality
\begin{eqnarray}
\eta_{m}^{2 m}+\frac{\eta_{m}}{1-\eta_{m}}+\frac{\eta_{m}^{m}}{\sqrt{1-\eta_{m}^{2 m}}}\leq 1.
\label{liu32}
\end{eqnarray}
\elem

\bpf
We first prove the uniqueness of the solution in $(0, 1/2)$ of the equation (\ref{liu24}).

Let $h(r)=r^{m}(2-3r) +2 r-1$. Then it is easy to note that $h(0)=-1<0$, $h(1/2)=(\frac{1}{2})^{m+1}>0$ and, for $r\in [0, 1/2]$, we have
$$h'(r) =m r^{m-1}(2-3 r)+2-3 r^{m}>0,
$$
showing that $h(r)$ is an increasing function of $r$ in $[0, 1/2]$, and thus, $h(r)=0$ has a unique root $\eta_{m}$ in $(0, 1/2)$.

Now we verify the inequality (\ref{liu32}). Let $y=\eta_{m}^{m}$. Then according to (\ref{liu24}), we have $\eta_{m}=(1-2y)/(2-3y)$.
Using this change of variables, we can rewrite (\ref{liu32}) in the following equivalent form:
\begin{eqnarray*}
y^{2}+\frac{1-2y}{1-y}+\frac{y}{\sqrt{1-y^{2}}}\leq 1 ~\mbox{ for }~ y\in[0,1/2).
\end{eqnarray*}
We note that
\begin{eqnarray*}
y^{2}+\frac{1-2y}{1-y}+\frac{y}{\sqrt{1-y^{2}}}-1=\frac{y[(-y^{2}+y-1)\sqrt{1-y^{2}}+1-y]}{(1-y)\sqrt{1-y^{2}}},
\end{eqnarray*}
and therefore, the inequality \eqref{liu32} is valid if and only if $(-y^{2}+y-1)\sqrt{1-y^{2}}+1-y\leq 0$ for $y\in[0,1/2)$, which holds
if and only if $y^{2}(y^{3}-y^{2}+y+1)\geq 0 $ for $y\in[0,1/2)$. Since the last inequality is obviously true, the proof is complete.
\epf
\section{Bohr-type inequalities for bounded analytic functions}\label{HLP-sec4}

\subsection{Proof of Theorem \ref{HLP-th1}}
Firstly, we consider the first part. Suppose that $f \in\mathcal{B}$,  $a:=|a_0|$ and $\omega \in\mathcal{B}_m$. Then, by
the classical Schwarz lemma and the Schwarz-Pick lemma or Lemma \ref{HLP-lem5}, we have
\begin{eqnarray}
|\omega(z)|&\leq & |z|^m, \quad z\in \mathbb{D},\label{liu41}\\
|f(u)|&\leq & \frac{|u|+a}{1+a|u|},\quad u\in \mathbb{D},\label{liu42}
\end{eqnarray}
which implies
\begin{eqnarray}
|f(\omega(z))|\leq  \frac{|\omega(z)|+a}{1+a|\omega(z)|}\leq \frac{r^{m}+a}{1+a r^{m}},\quad |z|=r<1.
\label{liu43}
\end{eqnarray}

According to Lemma \ref{HLP-lem2} 
and (\ref{liu43}), we have 
$$
B_1(f,r)+ \frac{1+ar}{(1+a)(1-r)} \|f_0\|_r^2 \leq \frac{(1-a^{2})r}{1-r},
$$
and thus,
\beqq
A_f(z) &\leq& 
 1-\left [1- \frac{r^{m}+a}{1+a r^{m}} -\frac{(1-a^{2})r}{1-r}\right ]= 1-\frac{(1-a)A_{m}(a, r)}{\left(1+a r^{m}\right)(1-r)},
\eeqq
where $A_m(a,r)=(1-r)(1-r^m)-(1+a)r(1+ar^m)$, which is clearly a decreasing function of $a\in [0,1]$.  Thus,
$$
A_m(a,r)\geq A_m(1,r)=(1-r)(1-r^m)-2r(1+r^m),
$$
and obtain that $A_f(z)\leq 1$ if $A_{m}(1, r)\geq 0$, which holds for $r\leq \alpha_{m}$, where $\alpha_{m}$ is the unique positive root in $(0, 1)$ of the equation
$A_{m}(1, r)=0$.

To show that the  radius $\alpha_{m}$ is best possible, we consider the functions
\be\label{HLP-eq1}
\omega(z)=z^m ~\mbox{ and }~\varphi_a(z)=\frac{z+a}{1+a z}=a+\left(1-a^{2}\right) \sum_{n=0}^{\infty}(-a)^{n} z^{n+1},~a \in[0,1).
\ee
For the two functions, we get that (for $z=r$)
\begin{eqnarray*}
A_{\varphi_a}(z)&=& |\varphi_a(z^m)| + B_1({\varphi_a},r)+ \frac{1+ar}{(1+a)(1-r)} \|{\varphi_a}_0-a\|_r^2\\
&=&\frac{a+r^{m}}{1+a r^{m}}+\frac{\left(1-a^{2}\right) r}{1-a r}+ \frac{\left(1-a^{2}\right)^{2} r^{2}}{(1+a)(1-r)(1-a r)}\\
&=&\frac{r^{m}+a}{1+a r^{m}}+\frac{r\left(1-a^{2}\right)}{1-r},
\end{eqnarray*}
and this expression is bigger than $1$ provided $r>\alpha_{m,a}$, where $\alpha_{m,a}$ is the unique positive root in $(0,1)$ of the equation
$A_{m,a}(r)=0$. Allowing $a\rightarrow 1^{-}$ gives that $\alpha_{m,1}=\alpha_{m}$ is the best possible constant.

Next, we prove the second part. As in the previous case, by Lemma \ref{HLP-lem2} and (\ref{liu43}), 
it follows easily that
$$B_f(z) \leq 1-\left [1-\left(\frac{r^{m}+a}{1+a r^{m}}\right)^{2}-\frac{(1-a^{2})r}{1-r}\right ]=1-\frac{\left(1-a^{2}\right) B_{m}(a, r)}{\left(1+a r^{m}\right)^{2}(1-r)},
$$
where $B_{m}(a, r)=(1-r)(1-r^{2m})-r(1+ar^m)^2$, 
which is clearly a decreasing function of $a\in [0,1]$.  Thus,
$$B_m(a,r)\geq B_m(1,r)=(1-r)(1-r^{2m})-r(1+r^m)^2= (1+r^m)(1-2r-r^m)=:B_m(r).
$$
We see that $B_f(z)\leq 1$ if $B_{m}(r)\geq 0$, which holds for $r\leq \beta_{m}$, where $\beta_{m}$ is the unique positive root in $(0, 1)$ of the equation $B_{m}(r)=0$, namely, \eqref{liu22} given by $1-2r-r^m=0$.

To show the radius $\beta_{m}$ is best possible, we consider the functions $\omega(z)=z^m$ and $\varphi_a(z)$ as above, and  find that (for $z=r$)
\beqq
B_{\varphi_a}(z)&= &|\varphi_a(z^m)|^2 + B_1({\varphi_a},r)+ \frac{1+ar}{(1+a)(1-r)} \|{\varphi_a}_0-a\|_r^2 \\
&=& \left(\frac{r^{m}+a}{1+a r^{m}}\right )^2+\frac{r\left(1-a^{2}\right)}{1-r},
\eeqq
and this expression is bigger than $1$ provided $r>\beta_{m,a}$, where $\beta_{m,a}$ is the unique positive root in $(0, 1)$ of the equation
$B_{m,a}(r)=0$. Allowing $a\rightarrow 1^{-}$ gives that $\beta_{m,1}=\beta_{m}$ is the
best possible constant.
Thus the proof of Theorem \ref{HLP-th1} is complete.
\hfill $\Box$

\subsection{Proof of Theorem \ref{HLP-th2}}

We begin to recall from Lemma \ref{HLP-lem1} that (as $\omega \in\mathcal{B}_m$ so that $|\omega(z)|\leq r^m$ and $f(\omega (0))=a_0$ with $a=|a_0|$),
\be\label{KayPon8-eq7}
|f(\omega (z))-a_0|\leq \sum_{k=1}^\infty |a_k|r^{mk} \leq
\left \{ \begin{array}{lr} \ds r^m\frac{1-a^2}{1-r^ma} & \mbox{ for $a \ge r^m$}, \\[4mm]
\ds r^m\frac{\sqrt{1-a^2}}{\sqrt{1-r^{2m}}} &\mbox{ for $a < r^m$}.
\end{array} \right .
\ee

For the first part of the proof of the theorem, we first consider $a\geq r^{m}$. Then it follows from (\ref{liu41}),
the first inequality on the right of \eqref{KayPon8-eq7} and Lemma \ref{HLP-lem2} that
\begin{eqnarray*}
C_f(z) 
&\leq & a+\frac{r\left(1-a^{2}\right)}{1-r}+\frac{r^{m}\left(1-a^{2}\right)}{1-a r^{m}}
= 1-\frac{(1-a) C_{m}(a, r)}{\left(1-a r^{m}\right)(1-r)},
\end{eqnarray*}
where
$$
C_{m}(a, r)=r^{m+1} a^{2}-\left(r-3 r^{m+1}+2 r^{m}\right) a-\left(r^{m}-r^{m+1}+2 r-1\right).
$$

Now, for fixed the value of $r$ in the expression, we obtain
$$\frac{\partial C_{m}(a, r)}{\partial a}=r\left [2a r^{m} -\left(1-3 r^{m}+2 r^{m-1}\right)\right ] \leq \left .\frac{\partial C_{m}(a, r)}{\partial a} \right |_{a=1} \, =
r\left [5r^{m} -2 r^{m-1}-1\right ],
$$
which is non-positive for $r\leq \nu_{m} $, $\nu_{m}\geq 3/5 $, where $\nu_{m} $ is the unique root in $(0, 1)$ of the equation $5 r^{m}-2 r^{m-1}-1=0$.
So, $C_{m}(a, r)$ is a decreasing function of $a\in [r^{m}, 1]$ and thus, we have
$$
C_{m}(a, r) \geq C_{m}(1, r)=5 r^{m+1}-3 r^{m}-3 r+1=:C_{m}(r).
$$
Clearly, $C(z)\leq 1$ if $C_{m}(r)\geq 0$, which holds for $r\leq \zeta_{m}$, where $\zeta_{m}$ is the unique positive root in $(0, 1/3)$ of the equation $C_{m}(r)=0$ from Lemma \ref{HLP-lem3}.

If $a< r^{m}\leq \zeta_{m}^{m} $, then combining (\ref{liu41}), Lemma \ref{HLP-lem2} and
the second inequality on the right in \eqref{KayPon8-eq7}, we have
\begin{eqnarray*}
C_f(z)
&\leq & a+\frac{r\left(1-a^{2}\right)}{1-r}+\frac{r^{m} \sqrt{1-a^{2}}}{\sqrt{1-r^{2 m}}}=:C(a, r).
\end{eqnarray*}
It is easy to see that
$$\frac{\partial C(a, r)}{\partial r}=\frac{1-a^{2}}{(1-r)^{2}}+m r^{m-1} \sqrt{\frac{1-a^{2}}{\left(1-r^{2 m}\right)^{3}}}>0,
$$
showing that $C(a, r)$ is monotonically increasing with respect to $r\in\left[0, \zeta_{m}\right]$ for each fixed $a\in [0,1)$. Thus,
we have from Lemma \ref{HLP-lem3} that
\begin{eqnarray*}
C_f(z) &\leq& a+\frac{\left(1-a^{2}\right) \zeta_{m}}{1-\zeta_{m}}+\frac{\zeta_{m}^{m} \sqrt{1-a^{2}}}{\sqrt{1-\zeta_{m}^{2 m}}}\\
&\leq&\zeta_{m}^{m}+\frac{\zeta_{m}}{1-\zeta_{m}}+\frac{\zeta_{m}^{m}}{\sqrt{1-\zeta_{m}^{2 m}}}\leq 1.
\end{eqnarray*}

To show the sharpness of the radius $\zeta_{m}$, we consider the functions $\omega$ and $\varphi_a$ as in \eqref{HLP-eq1}, and
obtain as before that (by setting $z=r$ for the first term in the definition of $C_{\varphi_a}(z)$)
\begin{eqnarray*}
C_{\varphi_a}(z)
&=&a+\frac{r\left(1-a^{2}\right)}{1-r}+\frac{r^{m}\left(1-a^{2}\right)}{1-a r^{m}},
\end{eqnarray*}
and the last expression shows the radius $\zeta_{m}$ is optimal.

For the proof of the second part of the theorem, when  $a\geq r^{m}$, it follows from (\ref{liu41}), Lemma \ref{HLP-lem2} and
the first inequality on the right of \eqref{KayPon8-eq7} that
$$D_f(z) \leq a^{2}+\frac{r\left(1-a^{2}\right)}{1-r}+\frac{r^{m}\left(1-a^{2}\right)}{1-a r^{m}}
=1-\frac{\left(1-a^{2}\right) D_{m}(a, r)}{\left(1-a r^{m}\right)(1-r)},
$$
where
$$
D_{m}(a, r)=-a r^{m}(1-2r)-r^{m}(1-r)-2 r+1.
$$
For $r\leq 1/2 $, it is clear that $D_{m}(a, r)$ is a decreasing function of $a$,  $a\in [r^{m},1]$. Hence
$$
D_{m}(a, r) \leq D_{m}(1, r)=-r^{m}(2-3r)-2 r+1=:D_{m}(r).
$$

Obviously, $D_f(z)\leq 1$ if $D_{m}(r)\leq 0$, which holds for $r\leq \eta_{m}\leq 1/2 $, where $\eta_{m}$ is the unique positive root in $(0, 1/2]$
of the equation $D_{m}(r)=0$  from Lemma \ref{HLP-lem4}.

If $a<r^{m} \leq \eta_{m}^{m}$, as in the previous case, we have
$$D_f(z) \leq a^{2}+\frac{r\left(1-a^{2}\right)}{1-r}+\frac{r^{m} \sqrt{1-a^{2}}}{\sqrt{1-r^{2 m}}}=:D(a, r).
$$
Since $D(a, r)$ is  clearly monotonically increasing with respect to $r\in[0, \eta_{m}]$, for each fixed value of $a$ in the expression,
it follows from Lemma \ref{HLP-lem4} that
\begin{eqnarray*}
D_f(z) &\leq& a^{2}+\frac{\left(1-a^{2}\right) \eta_{m}}{1-\eta_{m}}+\frac{\eta_{m}^{m} \sqrt{1-a^{2}}}{\sqrt{1-\eta_{m}^{2 m}}}\\
&\leq& \eta_{m}^{2 m}+\frac{\eta_{m}}{1-\eta_{m}}+\frac{\eta_{m}^{m}}{\sqrt{1-\eta_{m}^{2 m}}}\leq1.
\end{eqnarray*}
The sharpness part follows similarly. Thus, the proof of Theorem \ref{HLP-th2} is complete.\hfill $\Box$

\vspace{6pt}



It is a simple exercise to see that for $0\leq x\leq x_0 ~(\leq 1)$, we have
\be\label{HLP-eq3a}
\Phi (x)=x+A(1-x^2)\leq \Phi (x_0) ~\mbox{ whenever $0\leq A\leq 1/2$,}
\ee
and similarly,
\be\label{HLP-eq3b}
\Psi (x)=x^2+A(1-x^2)\leq \Psi (x_0) ~\mbox{ whenever $0\leq A\leq 1$.}
\ee
These two inequalities will be used in the proofs of Theorems \ref{HLP-th3} and \ref{HLP-th4}.

\subsection{Proof of Theorem \ref{HLP-th3}}
Firstly, we consider the first part. In view of (\ref{liu41}), (\ref{liu43}), Schwarz-Pick lemma and  \eqref{HLP-eq2}, we have
\begin{eqnarray*}
E_f(z) &\leq& |f(w(z))|+ \frac{r^{m}}{1-r^{2 m}}\left(1-|f(w(z))|^{2}\right)+\frac{\left(1-a^{2}\right) r^{2}}{1-r}\\
&\leq& \frac{r^{m}+a}{1+a r^{m}}+\frac{r^{m}}{1-r^{2 m}}\left[1-\left(\frac{r^{m}+a}{1+a r^{m}}\right)^{2}\right]+\frac{\left(1-a^{2}\right) r^{2}}{1-r} ~\mbox{ (by \eqref{HLP-eq3a})}\\
&=&\frac{r^{m}+a}{1+a r^{m}}+\frac{r^{m}\left(1-a^{2}\right)}{\left(1+a r^{m}\right)^{2}}+\frac{\left(1-a^{2}\right) r^{2}}{1-r}\\
&=&1+\frac{(1-a) E_{m}(a, r)}{\left(1+a r^{m}\right)^{2}(1-r)},
\end{eqnarray*}
for $r\in [0,\mu_{m}]$, since $\frac{r^{m}}{1-r^{2 m}}\leq\frac{1}{2}$ for $r\in [0,\mu_{m}]$, where $\mu_{m}=\sqrt[m]{\sqrt{2}-1}$ is the unique root in (0,1) of the equation $r^{2 m}+2 r^{m}-1=0$, and
\begin{eqnarray*}
E_{m}(a, r)&=&r^{2 m+2} a^{3}+ r^{m+2}\left(r^{m}+2\right) a^{2}+\left [ r^{2 m} (1-r) +r^{2}(2 r^{m}+1)\right ] a\\
&&+2 r^{m}(1- r) +r^{2}+r-1.
\end{eqnarray*}
For each fixed $r\in [0,1]$,
it is clear that $E_{m}(a, r)$ is a monotonically increasing function of $a\in  [0,1)$ and thus,
$$
E_{m}(a, r) \leq E_{m}(1, r)= r^{m}(r^m+2)[2r^{2} -r +1]+2 r^{2}+r-1 =:E_{m}(r).
$$
Therefore, $E_f(z)\leq 1$ if $E_{m}(r)\leq 0$, which is valid for $r\leq \gamma_{m}$,
where $\gamma_{m}$ is the unique positive root in $(0, 1)$ of the equation $E_{m}(r)=0$,
and obviously $\gamma_{m}< \mu_{m}.$

To show that the radius $\gamma_{m}$ is optimal, as in the
proofs of the previous two theorems, we consider the functions $\omega$ and $\varphi_a$ as in \eqref{HLP-eq1}, and
set $z=r$ for the first term in the definition of $E_{\varphi_a}(z)$ and obtain that
\begin{eqnarray}\label{HLP-eq4}
E_{\varphi_a}(z)
&=&\frac{a+r^{m}}{1+a r^{m}}+\frac{\left(1-a^{2}\right) r^{m}}{\left(1+a r^{m}\right)^{2}}+\frac{r^{2}\left(1-a^{2}\right)}{1-r} \nonumber\\
&=&1+\frac{(1-a)E_{m}(a,r)}{\left(1+a r^{m}\right)^{2}(1-r)},
\end{eqnarray}
which is larger than $1$ if and only if $E_m(a,r)>0$. Also, the expression on the right is smaller than or equal to $1$
for all $a\in [0,1)$, only in the case when $r\leq \gamma_{m}$. Finally, it also suggests that $a\to 1$ in \eqref{HLP-eq4}
shows that the expression on the right is larger than $1$ if $r>\gamma_{m}$.


Concerning the second sum in the theorem, it follows from (\ref{liu41}), (\ref{liu43}),
Schwarz-Pick lemma and \eqref{HLP-eq2} that
\begin{eqnarray*}
F_f(z)&\leq&  |f(w(z))|^2+ \frac{r^{m}}{1-r^{2 m}}\left(1-|f(w(z))|^{2}\right)+\frac{\left(1-a^{2}\right) r^{2}}{1-r}\\
&\leq& \left(\frac{r^{m}+a}{1+a r^{m}}\right)^{2}+\frac{r^{m}}{1-r^{2 m}}\left[1-\left(\frac{r^{m}+a}{1+a r^{m}}\right)^{2}\right]+\frac{\left(1-a^{2}\right) r^{2}}{1-r}
~\mbox{ (by \eqref{HLP-eq3b})}\\
&=&\left(\frac{r^{m}+a}{1+a r^{m}}\right)^{2}+\frac{r^{m}\left(1-a^{2}\right)}{\left(1+a r^{m}\right)^{2}}+\frac{\left(1-a^{2}\right) r^{2}}{1-r}\\
&=&1+\frac{\left(1-a^{2}\right) F_{m}(a, r)}{\left(1+a r^{m}\right)^{2}(1-r)},
\end{eqnarray*}
for $r \in [0,\tau_{m}]$, since $\frac{r^{m}}{1-r^{2 m}} \leq 1$ for $r \in [0,\tau_{m}]$, where $\tau_{m}=\sqrt[m]{\frac{\sqrt{5}-1}{2}}\,$ is the unique root
in $(0, 1)$ of the equation $r^{2m}+r^{m}-1=0 $, and
$$F_{m}(a, r)=(r^{2m}+r^{m}-1)(1-r) +r^2(1+a r^{m})^2.
$$
Since $F_{m}(a, r)$ is clearly an increasing function of $a$ in $[0,1)$, it follows that
$$F_{m}(a, r) \leq F_{m}(1, r)=r^{m}(r^m+1)[r^{m} -r +2]+ r-1  
=:F_{m}(r).
$$
Thus, $F_f(z)\leq 1$ if $F_{m}(r)\leq 0$, which holds for $r\leq \delta_{m}$, where $\delta_{m}$ is the unique positive root in $(0, 1)$ of the
equation $F_{m}(r)=0$, and obviously, $\delta_{m}<\tau_{m}$. Sharpness part may be proved similarly. Thus, we conclude the proof of the theorem.
\hfill $\Box$

\subsection{Proof of Theorem \ref{HLP-th4}}
Firstly, we consider the first part. Clearly, $\frac{2r}{1-r^{2 m}} \leq 1$ if  $r\in [0,\xi_{m}]$, where $\xi_{m}$ is the unique root in $(0,1)$
of the equation $r^{2 m}+2 r-1=0$. As before, it follows from (\ref{liu41}), (\ref{liu43}), the Schwarz-Pick lemma and \eqref{HLP-eq2} that
\begin{eqnarray*}
G_f(z) &\leq& |f(w(z))|+\frac{r}{1-r^{2 m}}\left(1-|f(w(z))|^{2}\right)+\frac{\left(1-a^{2}\right) r^{2}}{1-r}\\
&\leq& \frac{r^{m}+a}{1+a r^{m}}+\frac{r}{1-r^{2 m}}\left[1-\left(\frac{r^{m}+a}{1+a r^{m}}\right)^{2}\right]+\frac{\left(1-a^{2}\right) r^{2}}{1-r} ~\mbox{ (by \eqref{HLP-eq3a})}\\
&=&\frac{r^{m}+a}{1+a r^{m}}+\frac{r\left(1-a^{2}\right)}{\left(1+a r^{m}\right)^{2}}+\frac{\left(1-a^{2}\right) r^{2}}{1-r}\\
&=&1+\frac{(1-a) G_{m}(a, r)}{\left(1+a r^{m}\right)^{2}(1-r)},
\end{eqnarray*}
for $r\in [0,\xi_{m}]$, where
\begin{eqnarray*}
G_{m}(a, r)&=&r^{2 m+2} a^{3}+r^{m+2}(r^{m}+2) a^{2}+ [r^{2m}(1-r)+ r^{m+1}(2r+1)+r(1-r^{m-1})] a\\
&&+r^{m}- r^{m+1}+2r-1,
\end{eqnarray*}
which, for each fixed value of $r\in[0,1]$,
is clearly monotonically increasing with respect to $a\in  [0,1)$, because the coefficients of $a^3,\,a^2$ and $a$ are non-negative for $r\in[0,1]$. Hence
$$
G_{m}(a, r) \leq G_{m}(1, r)=2 r^{2 m+2}-r^{2 m+1}+r^{2 m}+4 r^{m+2}+3r-1=:G_{m}(r).
$$
Therefore, $G_f(z)\leq 1$ if $G_{m}(r)\leq 0$, which is valid for $r\leq \theta_{m}$, where $\theta_{m}$ is the unique positive root in $(0, 1)$
 of the equation $G_{m}(r)=0$, and it is clear that  $\theta_{m}<\xi_{m}$.

To show that the radius $\theta_{m}$ is optimal,  we consider the functions $\omega$ and $\varphi_a$ as in \eqref{HLP-eq1}, and
set $z=r$ for the first term in the definition of $G_{\varphi_a}(z)$ and obtain that
\begin{eqnarray*}
G_{\varphi_a}(z))
&=&\frac{a+r^{m}}{1+a r^{m}}+\frac{\left(1-a^{2}\right) r}{\left(1+a r^{m}\right)^{2}}+\frac{r^{2}\left(1-a^{2}\right)}{1-r}\\
&=& 1+\frac{(1-a)G_{m}(a,r)}{\left(1+a r^{m}\right)^{2}(1-r)},
\end{eqnarray*}
which is larger than $1$ if $G_{m}(r)> 0$, and using the earlier arguments, this is valid for $r> \theta_{m}$,
where $\theta_{m}$ is the unique positive root in $(0, 1)$  of the equation $G_{m}(r)=0$.

Obviously, $1-\frac{r}{1-r^{2 m}} \geq 0$ if $r \in [0,\chi_{m}]$, where $\chi_{m}$ is the unique root in (0,1) of the equation $r^{2m}+r-1=0 $.
 Again, as in the previous case, we have for $r \in [0,\chi_{m}]$ that
\begin{eqnarray*}
H(z)&\leq& |f(w(z))|^{2}+\frac{r}{1-r^{2 m}}\left(1-|f(w(z))|^{2}\right)+\frac{\left(1-a^{2}\right) r^{2}}{1-r}\\
&\leq& \left (\frac{r^{m}+a}{1+a r^{m}}\right )^{2}+\frac{r}{1-r^{2 m}}\left[1-\left(\frac{r^{m}+a}{1+a r^{m}}\right)^{2}\right]+\frac{\left(1-a^{2}\right) r^{2}}{1-r} ~\mbox{ (by \eqref{HLP-eq3b})}\\
&\leq& \left(\frac{r^{m}+a}{1+a r^{m}}\right)^{2}+\frac{r\left(1-a^{2}\right)}{\left(1+a r^{m}\right)^{2}}+\frac{\left(1-a^{2}\right) r^{2}}{1-r}\\
&=&1+\frac{\left(1-a^{2}\right) H_{m}(a, r)}{\left(1+a r^{m}\right)^{2}(1-r)},
\end{eqnarray*}
for $r \in [0,\chi_{m}]$, where
$$
H_{m}(a, r) =(1-r)(r^m+r-1) +r^2(1+a r^{m})^{2}.
$$
As $H_{m}(a, r)$ is an increasing function for $a$ in $ [0,1)$, we obtain that
$$H_{m}(a, r) \leq H_{m}(1, r)=r^{2 m+2}-r^{2 m+1}+r^{2 m}+2 r^{m+2}+2r-1=:H_{m}(r).
$$
Thus, $H(z)\leq 1$ is valid if $H_{m}(r)\leq 0$, which holds for $r\leq \vartheta_{m}$, where $\vartheta_{m}$ is the unique positive root in $(0, 1)$ of the equation $H_{m}(r)=0$,
and $\vartheta_{m}<\chi_{m}$ holds obviously. The sharpness part is similar. The proof of the theorem is complete. \hfill $\Box$

\subsection{Proof of Theorem \ref{HLP-th5}} By Theorem \Ref{Theo-F} (and its proof), we only need to prove the case $m\geq 2$.

In fact, if $m\geq 2$ and $ a \geq r^{m}$, then it follows from (\ref{liu41}), Lemmas \ref{HLP-lem1} and \ref{HLP-lem2} that
\begin{eqnarray*}
I_{f}(z) &\leq & a+\frac{r\left(1-a^{2}\right)}{1-r}+{\left[\frac{|\omega(z)| (1-a^{2})}{1-a|\omega(z)|}\right]^{2}}
\leq a+\frac{r\left(1-a^{2}\right)}{1-r}+\frac{r^{2 m}\left(1-a^{2}\right)^{2}}{\left(1-a r^{m}\right)^{2}}\\
&=&1+\frac{(1-a)I_{m}(a,r)}{(1-r)(1-a r^{m})^{2}},
\end{eqnarray*}
where
\begin{eqnarray*}
I_{m}(a,r)&=& (2r-1+ar)(1-2ar^m+a^2 r^{2m}) +r^{2 m}(1-r)(1-a^{2})(1+a)\\
&=&  -r^{3 m}(1-2r)a^{3}-[r^{2 m}(2-3 r)+2 r^{m+1}]a^{2}\\
&&+[r^{2m}(1-r)+2(1-2r)r^{m}+r]a  + (1-r)r^{2m}+2r-1.
\end{eqnarray*}
As
$$\frac{\partial^{2} I_{m}(a, r)}{\partial a^{2}} =-6r^{3 m}(1-2r)a^{3}-2[r^{2 m}(2-3 r)+2 r^{m+1}] \leq 0 ~\mbox{ for $r\leq 1/2$,}
$$
it follows that for $ a \geq r^{m}$
$$ \frac{\partial^{2} I_{m}(a, r)}{\partial a^{2}} \leq\left.\frac{\partial^{2} I_{m}(a, r)}{\partial a^{2}}\right|_{a=r^{m}}
=-6 r^{3 m}(1-2 r)-2 r^{2 m}(2-3 r)-4 r^{m+1}\leq 0,
$$
showing that $\frac{\partial I_{m}(a, r)}{\partial a}$ is a decreasing function of $a\in [r^{m},1]$. Thus, we obtain
$$
\frac{\partial I_{m}(a, r)}{\partial a} \geq\left.\frac{\partial I_{m}(a, r)}{\partial a}\right|_{a=1}=r\left(11 r^{2 m}-6 r^{2 m-1}-8 r^{m}+2 r^{m-1}+1\right),
$$
which is non-negative for $r\leq \xi_{m}$, where $\xi_{m}\approx 0.487478 $ is the unique root in $(0,1)$ of the equation $11r^{2m}-6r^{2m-1}-8r^{m}+2r^{m-1}+1=0$.
Therefore, $I_{m}(a,r)$ is clearly monotonically increasing  with respect to $a\in [r^{m}, 1]$, and thus, we see that
$$
I_{m}(a, r) \leq I_{m}(1, r)=(3 r-1)\left(1-r^{m}\right)^{2} \leq 0 ~\mbox{ for $r\leq 1/3$.}
$$
We conclude that $I_{f}(z)\leq 1$ and hence,  (\ref{liu29}) holds,  for $r\leq 1/3$ and $ a \geq r^{m}$.

Next, we observe that $I_{m}(a,r)$ is an increasing function of $r\in [0,1)$ and may be written as
\begin{eqnarray*}
I_{m}(a,r) &=& 1+\frac{1-a}{(1-r)(1-a r^{m})^{2}}[(2 a^{3}+3 a^{2}-a-1) r^{2 m+1}-(a^{3}+2 a^{2}-a-1) r^{2 m}\\
&&-(2 a^{2}+4 a) r^{m+1}+2 a r^{m}+(a+2) r-1].
\end{eqnarray*}
Now, for  $r\leq 1/3 $, we have
\begin{eqnarray*}
I_f(z)
&\leq& I_{m}(a,1/3)=1+\frac{1-a}{2\left(3^{m}-a\right)^{2}}\left[-a^{3}-3 a^{2}+2 a+2+3^{m}(1-a)\left(2 a-3^{m}\right)\right]\\
&\leq& 1+\frac{1-a}{2\left(3^{m}-a\right)^{2}}\left[-a^{3}-3 a^{2}+2 a+2+3(1-a)(2 a-3)\right]\\
&=&1+\frac{(1-a)^2}{2\left(3^{m}-a\right)^{2}}\left(a^{2}+10a-7\right)\leq 1,
\end{eqnarray*}
since $a^{2}+10a-7\leq 0$ if $ 0\leq a\leq a^{*}=-5 +4\sqrt{2}$. In the third inequality above we have used the fact that
$3^{m}(2 a-3^{m})\leq 3(2 a-3)$, i.e. $2a-3^m-3\leq 0$ for all $m\geq 1$. 

Finally, if $m\geq 2$ and $0 \leq a < r^{m}\leq (\frac{1}{3})^{m}$, then it follows from (\ref{liu41}), Lemmas \ref{HLP-lem1} and \ref{HLP-lem2} that
\begin{eqnarray*}
I_{f}(z) &\leq& a+\frac{r\left(1-a^{2}\right)}{1-r}+{\left[\frac{|\omega(z)| \sqrt{1-a^{2}}}{\sqrt{1-|\omega(z)|^{2}}}\right]^{2}}\\
&\leq & a+\frac{r\left(1-a^{2}\right)}{1-r}+\frac{r^{2 m}\left(1-a^{2}\right)}{1-r^{2 m}}=:J_{m}(a,r).
\end{eqnarray*}
We notice that $J_{m}(a,r)$ is monotonically increasing with respect to $r$ and thus, for $0 \leq a < r^{m}\leq (\frac{1}{3})^{m}$, we obtain that
\begin{eqnarray*}
I_{f}(z)&\leq& J_{m}(a,1/3)=a+\frac{1}{2}\left(1-a^{2}\right)+\frac{\left(\frac{1}{3}\right)^{2 m}\left(1-a^{2}\right)}{1-\left(\frac{1}{3}\right)^{2 m}}\\
&\leq&\left(\frac{1}{3}\right)^{m}+\frac{1}{2}+\frac{\left(\frac{1}{3}\right)^{2 m}}{1-\left(\frac{1}{3}\right)^{2 m}}
\\ &\leq&
\frac{1}{9}+\frac{1}{2}+\frac{1}{80}<1.
\end{eqnarray*}

To show that the radius  is optimal,  we consider the functions $\omega$ and $\varphi_a$ as in \eqref{HLP-eq1}, and
set $z=r$ for the first term in the definition of $I_{\varphi_a}(z)$ and obtain that
\begin{eqnarray*}
I_{\varphi_a}(z)
&=&a+\frac{r\left(1-a^{2}\right)}{1-r}+\frac{r^{2 m}\left(1-a^{2}\right)^{2}}{\left(1-a r^{m}\right)^{2}},
\end{eqnarray*}
and the last expression easily delivers a proof of the sharpness part. The proof of Theorem \ref{HLP-th5} is complete.\hfill $\Box$

\subsection*{Acknowledgments}
This research of the first two authors are partly supported by Guangdong Natural Science Foundations (Grant No. 2018A030313508).
The work of the third author was supported by Mathematical Research Impact Centric Support (MATRICS) of
the Department of Science and Technology (DST), India  (MTR/2017/000367). The authors of this paper thank the referees very much for their valuable comments and suggestions to this paper.

%

\end{document}